\documentclass[10pt]{article}
\usepackage{array}
\usepackage[mathscr]{euscript}
\usepackage{amsmath,amssymb,amsbsy,amsfonts,amsthm,latexsym,verbatim}
\usepackage{graphicx,color,dsfont}
\usepackage{epsfig, tkz-berge,tikz,pgf,tkz-graph}
\usepackage{float}
\usepackage{algorithm}
\usepackage{color}
\usepackage[noend]{algpseudocode}
\usepackage{centernot}
\usepackage{caption}
\usepackage{multicol}
\usepackage{enumerate}
\usepackage{relsize}
\usepackage[T1]{fontenc}
\usepackage[utf8]{inputenc}
\usepackage[numbers, sort&compress]{natbib}
\usepackage{longtable}
\usepackage{enumitem}

\definecolor{burgundy}{rgb}{0.65,0.1,0.43}
\definecolor{purple}{rgb}{0.25,0.1,0.63}
\definecolor{fresia}{rgb}{0.85,0.01,0.83}
\definecolor{lightblue}{rgb}{0.05,0.81,0.83}
\definecolor{gold}{rgb}{0.85,0.51,0.23}

\def\qed{~\hfill\vrule height .9ex width .8ex depth -.1ex}

\usetikzlibrary{shapes,shapes.geometric,arrows,fit,calc,positioning,automata,}

\newcommand{\pf}{\medskip\noindent{\bf Proof:}\hspace{.7em}}

\newcommand{\beq}{\begin{equation}}
\newcommand{\eeq}{\end{equation}}
\newcommand{\bea}{\begin{eqnarray*}}
\newcommand{\eea}{\end{eqnarray*}}
\newcommand{\diam}{\mathop{\rm diam}}
\newcommand{\rad}{\mathop{\rm rad}}
\newtheorem{theorem}{Theorem}

\newtheorem{lemma}{Lemma}
\newtheorem{remark}{Remark}
\newtheorem{definition}{Definition}
\newtheorem{proposition}{Proposition}

\begin{document}
\sloppy

\title{The Neighbor Matrix: \\Generalizing A Graph's Degree Sequence}
\author{
Jonathan W. Roginski$^*$ \\Department of Applied Mathematics\\ Naval Postgraduate School, Monterey, CA \\ $^*$Corresponding author: jonathan.w.roginski.mil@mail.mil\bigskip \\Ralucca M. Gera\\Department of Applied Mathematics \\ Naval Postgraduate School, Monterey, CA\bigskip \\ Erik C. Rye\\ Department of Applied Mathematics\\ Naval Postgraduate School, Monterey, CA}

\date{}

\maketitle

\begin{abstract}
{The newly introduced neighborhood matrix extends the power of adjacency and distance matrices to describe the topology of graphs. The adjacency matrix enumerates which pairs of vertices share an edge and it may be summarized by the degree sequence, a list of the adjacency matrix row sums. The distance matrix shows more information, namely the length of shortest paths between vertex pairs. We introduce and explore the neighborhood matrix, which we have found to be an analog to the distance matrix what the degree sequence is to the adjacency matrix. The neighbor matrix includes the degree sequence as its first column and the sequence of all other distances in the graph up to the graph's diameter, enumerating the number of neighbors each vertex has at every distance present in the graph. We prove this matrix to contain eleven oft-used graph statistics and topological descriptors. We also provide insight into two applications that show potential utility of the neighbor matrix in comparing graphs and identifying topologically significant vertices in a graph.}
\bigskip\\
2010 Math Subject Classification: 05C07, 05C12, 05C50
\end{abstract}

\section{Introduction}  \label{section:motivation}
Matrices such as the adjacency, distance, reciprocal distance, walk, reachability, and Laplacian offer well-studied, compact structures that represent graph information \citep{Harary1969, CZ, newman, EDM1, EDM2, mihalic1992distance, RDM, axelrod2015structure, merris1994laplacian}. Two of these matrix representations that form the foundation of this work are the adjacency matrix and the distance matrix. The adjacency matrix enumerates which pairs of vertices share an edge. The degree sequence is one summary of the adjacency matrix: a list of the adjacency matrix row sums. The distance matrix shows more information, namely the length of a shortest path between vertex pairs. We introduce and explore the neighbor matrix, a new matrix that extends the degree sequence to capture the distribution of each distance in a graph from $1$ to $k$, where $k$ is the graph diameter.  

\begin{figure}[ht]

\begin{multicols}{2}
    \centerline{\includegraphics[width=0.45\textwidth]{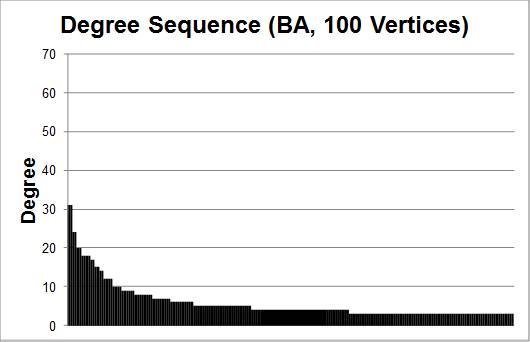}}
    
    \columnbreak

    
    \centerline{\includegraphics[width=0.45\textwidth]{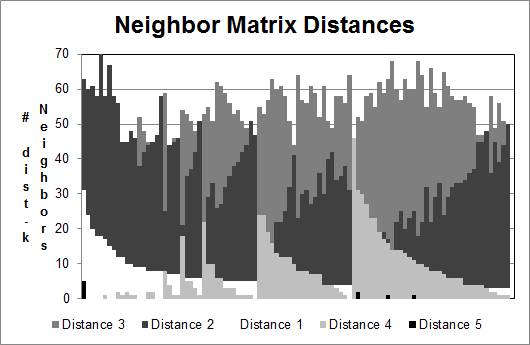}}

\end{multicols}

\caption{Comparing Topological Information of Neighbor Matrix to Degree Sequence}
\label{fig:NeighborVDegree}
\end{figure}

We need only consider Figure~\ref{fig:NeighborVDegree} to see the topological richness of the neighbor matrix. On the left in the figure is a depiction of the degree sequence of a Barab{\'a}si-Albert graph on $100$ vertices. The right figure depicts the distances accounted for by the neighbor matrix. The $distance-1$ neighbors are exactly the degree sequence of the graph, depicted in white on the right. The illustration shows $distance-2$ neighbors in the background of the right picture in medium gray, followed by $distance-3$ neighbors in dark gray, and $distance-4$ neighbors in light gray. We depict the small number of $distance-5$ neighbors in the foreground of the picture, in black. 

We use combinatoric techniques to show the neighbor matrix to include many of the statistics and topological characteristics currently used to describe graphs. Finally, we describe the use of the neighbor matrix in two applications that highlight its potential usefulness. We show that where current topological descriptors fail to discriminate between non-isomorphic graphs, the neighbor matrix may enable a graph comparison mechanism. We also present a methodology that leverages the topological information resident in the neighbors matrix to improve upon analysis using currently-defined centrality measures that seek to identify those vertices in the graph which exert the most influence over a graph's structure.


\section{Definitions}  \label{section:definitions}
We first review established definitions of graph theory and complex network terms that are fundamental to the work described in this paper. 

This work is limited to analysis of {\it simple, connected} graphs (no multiple edges nor loops) $G = (V(G), E(G))$, where $V(G)$ is the vertex set and $E(G)$ the edge set. We annotate the vertices of a graph $v_1, v_2, \dots, v_n$, where $n$ is the number of vertices in the graph. Two vertices are called \textit{adjacent} if there is an edge between them, otherwise they are \textit{nonadjacent}. The edge set $E(G)$ contains each adjacency between vertices. A {\textit cycle} is a sequence of adjacent vertices that begins and ends at the same vertex, but repeats no other vertex. A graph with no cycles is called {\textit acyclic} and all connected, acyclic graphs are called {\textit trees}. A graph $H$ is a {\it subgraph} of $G$ (denoted $H \subseteq G$) if $V(H) \subseteq V(G)$ and $E(H) \subseteq E(G)$. If $H \subseteq G$ and either $V(H)$ is a proper subset of $V(G)$ or $E(H)$ is a proper subset of $E(G)$, then $H$ is a proper subgraph of $G$, $ H \subset G$.  
    
Two graphs, $G$ and $H$ are {\it isomorphic} if there is a one-to-one correspondence $\phi$ from $V(G)$ to $V(H)$ such that $uv \in E(G)$ if and only if $\phi(u)\phi(v) \in E(H)$. A {\it graph invariant} is a property that has the same value for every pair of isomorphic graphs. An {\it automorphism} is an isomorphism from a graph $G$ to itself. Since the identity and inverse are both automorphisms and the composition of two automorphisms is itself an automorphism, the set of all automorphisms of a graph $G$ forms a {\it group} under the operation of composition. We denote the {\it automorphism group} of $G$ to be $Aut(G)$. Suppose $v$ is a vertex of graph $G$. The set of all vertices to which $v$ may be mapped by an automorphism of $G$ in an {\it orbit} of $G$. The automorphism relating two vertices is an equivalence relation resulting in equivalence classes that are the {\it orbits} of $G$. We define the orbit of vertex $i$ in graph $G$ as $o(i_G)$. If $G$ contains a single orbit, then $G$ is {\it vertex-transitive}.

The \textit{adjacency matrix} $A$ is comprised of entries $a_{ij}, \forall i,j \leq n$, where $n$ is the number of vertices in the graph. An entry $a_{ij} = 1$ represents the adjacency of vertex $i$ and vertex $j$; $a_{ij} = 0$ otherwise. 
A shortest path between two nonadjacent vertices is called a \textit{geodesic}. The length of a geodesic between two vertices is the \textit{distance} between the vertices. The \textit{distance matrix} $D$ is comprised of entries $d_{ij}, \forall i,j \leq n$, where $n$ is the number of vertices in the graph and $d_{ij}$ represents distance between vertex $i$ and vertex $j$. The \textit{average distance in $G$} is computed as the fraction of all pairwise distances out of all possible distances, $\displaystyle \frac{2}{n(n-1)} \sum_{i\neq j} d(v_i, v_j).$ By convention, $d(x,y) = 0$ if $x$ and $y$ are in different components. For a given vertex $v$ in graph $G$, $u$ is a \textit{$k$-step (or $k$-hop) neighbor} of $v$ if $d(u, v) \leq k$ for $1 \le k \le n-1$. Wu and Dai introduced a specification of this measure in \citep{exact_k}: a vertex $u$ is an \textit{exact $k$-hop neighbor} of $v$ if $d(u, v) = k$ for $1 \le k \le n-1$. 

Let $v$ be a vertex of $G$. The \textit{eccentricity} of the vertex $v$, $e(v)$, is the maximum distance from $v$ to any other vertex in $G$: $e(v)=\max\{d(v,w) : w \in V(G)\}.$
The \textit{radius of $G$} is the minimum eccentricity among the vertices of $G$: $\rad(G)=\min\{e(v) : v \in V(G)\}.$
The \textit{diameter of $G$} is the maximum eccentricity among the vertices of $G$: $\diam(G)=\max\{e(v) : v \in V(G)\}.$
The \textit{center of $G$}, $Cen(G)$, is the subgraph induced by those vertices of $G$ having minimum eccentricity:
$\displaystyle{Cen(G)=G[\{v \in V(G): e(v)=\rad(G)\}]}.$
The \textit{periphery of $G$}, $Per(G)$, is the subgraph induced by those vertices of $G$ having maximum eccentricity:
$Per(G)=G[\{v \in V(G): e(v)=\diam(G)\}].$
The \textit{$k^{th}$ power} of an undirected graph $G$ is the graph $G^k  = (V(G), E(G^k))$, where $E(G^k) = \{uv : d_G(u, v) \le k\}.$

A vertex $v_i$'s \textit{degree} denotes the number of vertices to which $v_i$ is adjacent. The \textit{degree sequence of $G$} of a graph is an integer sequence $d_1, d_2, \dots, d_n$ where $n = |V(G)|$ and $d_i$ is the degree of vertex $i$. We will use the convention that the sequence is non-increasing. The \textit{degree distribution of $G$} is the probability distribution of the degrees of the nodes in the graph, i.e. what fraction of the vertices have degree $k$ ($\delta(G) \le k \le \Delta(G)$), where $\delta(G)$ and $\Delta(G)$ are the \textit{minimum and maximum degree} in $G$, respectively.
The \textit{average degree of $G$} is 
$\displaystyle{\frac{\sum_{v_i \in V(G)} deg(v_i)}{n} =}$
$\displaystyle{\frac{2m}{n}}$, where $m$ is the number of edges in the graph. The \textit{density} of a graph is the ratio of possible edges to the edges that are actually present in the graph: 
$\displaystyle{\frac{m}{\binom{n}{2}}  = \frac{2m}{n(n-1)}}.$

The {\it average clustering coefficient} is the ratio of triangles in the graph to the number of connected triples (i.e., connected subgraphs on three vertices): $$ \text{Average Clustering Coefficient = } \frac{\text{number of triangles x 3}}{\text{number of connected triples}}.$$  {\it Pearson's correlation coefficient} in graphs is a measure of {\it assortative mixing}, i.e., the extent vertices with high degree are adjacent to each other. See \cite{newman2002assortative} for a detailed treatment of the topic. The {\it s-metric} determines the extent to which the graph being examined has a ``hub-like'' core; see \citep{smetric1,smetric2}.

A variety of ``centrality'' measures serve to provide insight into which vertices are the most \textit{influential} in a graph \citep{newman}. The word influential is emphasized because there is no commonly accepted definition in a graph topological context (and one is not offered here). We will leverage a basic understanding of ``importance'' in this paper: a vertex with high centrality is more ``important'' in some sense than a vertex of smaller centrality. The {\it degree} of vertex $i$ is sometimes referred to as {\it degree centrality} of vertex $i$ (degree centrality is often normalized through division by $n-1$). {\it Closeness centrality} is a measure of the distance from a vertex to all other vertices, calculated as the inverse mean distance. {\it Betweenness centrality} is the extent to which a given vertex lies on the shortest paths to other vertices.

The reader may refer to~\citep{Harary1969, CZ} for additional graph theory terminology and to~\citep{newman} for additional complex network terminology.

\bigskip 

We now introduce the following definitions as alternative (and more intuitive) terminology to the exact $k$-hop neighbor. The definitions will facilitate development of the neighbor matrix. 

\begin{definition} 
    Let $G$ be a graph with $u, v \in V(G)$. Vertex $u$ is a \emph{distance-$k$ neighbor} of $v$ if and only if $d(u, v) = k$, where $k = diam(G) (1 \le k \le n-1$). 
\end{definition}


\begin{definition}
    Let $G$ be a graph, with $v \in V(G)$. The \emph{distance-$k$ neighborhood} of $v$ is $$N^{dist-k}(v) = \{u \in V(G) : d(u,v) = k\}.$$
\end{definition}

\section{The Neighbor Matrix of Graph G: $X^{dist-k}(G)$} \label{section:matrix}

We introduce a matrix that is a multi-dimensional analog of the degree sequence. As the degree sequence summarizes the adjacencies in the adjacency matrix, the neighbor matrix summarizes the distributions of each vertex-vertex geodesic length found in the distance matrix.

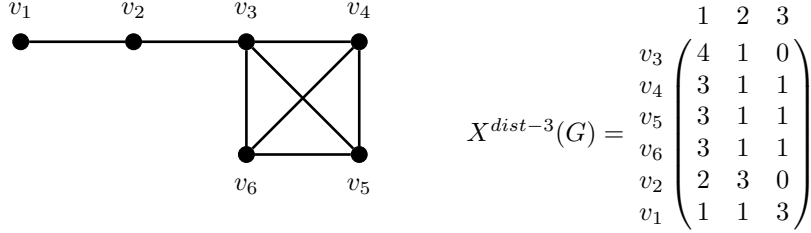
\begin{figure}[ht]
        \begin{multicols}{2}
        \centering
            \begin{tikzpicture}
            \tikzstyle{every node} = [inner sep = 0.0, minimum size = 6pt, draw = black]
            \node (y) [color = white,fill = black, circle , draw = black] at (0.00,19.00) [label= {[label distance = 2mm]$v_{1}$}] {};
            \node (x) [color = white,fill = black, circle , draw = black] at (1.50,19.00) [label= {[label distance = 2mm]$v_{2}$}] {};
            \node (u1) [color = white,fill = black, circle , draw = black] at (3.00,19.00) [label= {[label distance = 2mm]$v_{3}$}] {};
            \node (u2) [color = white,fill = black, circle , draw = black] at (4.50,19.00) [label= {[label distance = 2mm]$v_{4}$}] {};
            \node (u3) [color = white,fill = black, circle , draw = black] at (4.50,17.50) [label= {[label distance = 2mm]below:$v_{5}$}] {};
            \node (u4) [color = white,fill = black, circle , draw = black] at (3.00,17.50) [label= {[label distance = 2mm]below:$v_{6}$}] {};
            \draw [black,line width = 1]  (y) -- (x);
            \draw [black,line width = 1]  (x) -- (u1);
            \draw [black,line width = 1]  (u1) -- (u2);
            \draw [black,line width = 1]  (u1) -- (u3);
            \draw [black,line width = 1]  (u1) -- (u4);
            \draw [black,line width = 1]  (u2) -- (u3);
            \draw [black,line width = 1]  (u2) -- (u4);
            \draw [black,line width = 1]  (u3) -- (u4);
            \end{tikzpicture}
        \columnbreak
        \noindent
            $ X^{dist-3}(G) =
            \bordermatrix{
            ~ & 1 & 2 & 3 \cr
            v_3 & 4 & 1 & 0 \cr
            v_4 & 3 & 1 & 1 \cr
            v_5 & 3 & 1 & 1 \cr
            v_6 & 3 & 1 & 1 \cr 
            v_2 & 2 & 3 & 0 \cr
            v_1 & 1 & 1 & 3 \cr}$
            \vspace{.6cm}
        \end{multicols}
    \caption{Graph $G$ with associated neighbor matrix}
    \label{fig:example}
    \end{figure}

The neighbor matrix is a multi-dimensional analog of the degree sequence. As the degree sequence summarizes the adjacencies in the adjacency matrix, the neighbor matrix summarizes the distributions of each vertex-vertex geodesic length found in the distance matrix. See Definition~\ref{def:NeighborMatrix} and Figure~\ref{fig:example}.

\bigskip 
\begin{definition}\label{def:NeighborMatrix}
    Let $G$ be a graph with $V(G) = \{v_1, v_2, \ldots, v_n\}$. The neighbor matrix, $$X^{dist-k}(G) =  [x_{ij}], (1 \le i \le n, 1 \le j\le k),$$ where  $n = |V(G)|$, $k = diam(G)$  $(1 \le k \le n-1)$ and $x_{ij} = |N^{dist-j}(i)|$, $i \in V(G); \, 1\leq j \leq k.$ We sort the rows of $X^{dist-k}$ in a reverse lexicographic manner by organizing the rows of $X^{dist-k}(G)$ in non-increasing order of the entries in the first column. If there is a tie in some column $j<k$, we sort affected rows in non-increasing order of column $j+1$. 
\end{definition}

The entries $x_{ij}$ represent the count of dist-$j$ neighbors of vertex $i$. This is the same number defined by Bloom, et al, and later by Buckley and Harary as the $d_{ij}^{th}$ entry of the \textit{distance degree sequence} \cite{bloom1983some, buckley1990distance}. In addition, we define $X_i^{dist-k}(G)$ to be row $i$ of the neighbor matrix associated with graph $G$, the row associated with $v_i \in V(G)$.

    \setlength{\columnsep}{0pt}
    Figure~\ref{fig:example} depicts a graph $G$ and its corresponding neighbor matrix $X^{dist-3}(G)$.

\begin{lemma}
\label{lemma:uniqueness}
There is a unique neighbor matrix $X^{dist-k}(G)$ for each graph $G$.
\end{lemma}

\pf{
Suppose otherwise that a given graph can have two different neighbor matrices. Let $K^1 = X_1^{dist-k}(G)$ and $K^2 = X_2^{dist-k}(G)$ such that $K^1 \neq K^2$. Then, there is an entry $k_{ij}^1 \neq k_{ij}^2$. By construction and ordering of the neighbor matrix, vertex $i$ in $G$ simultaneously has two different numbers of neighbors at distance $j$, which is not possible.~\qed{}
}

Given Lemma~\ref{lemma:uniqueness}, the following observation is immediate:
\begin{proposition} 
\label{Prop:Congruent}
    Let $G$ and $H$ be two graphs. If $G \cong H$, then $X^{dist-k}(G) = X^{dist-k}(H).$
\end{proposition}

\pf{
By definition of isomorphism, there is a one-to-one correspondence $\phi$ from $V(G)$ to $V(H)$ such that $uv \in E(G)$ if and only if $\phi(u)\phi(v) \in E(H)$. Therefore, the distance in $G$ from vertex $i$ to vertex $j$ is the same as the distance in $H$ from $\phi(i)$ to $\phi(j), \ \ \forall i,j \in G$. It follows directly by construction of the neighbor matrix that $\displaystyle{X^{dist-k}(G) = X^{dist-k}(H)}.$~\qed{}
}

\begin{remark}
\label{rmk:different}
    The converse of Theorem$~\ref{Prop:Congruent}$ is false.  That is, $$X^{dist-k}(G) = X^{dist-k}(H) \centernot\implies G \cong H.$$ 
\end{remark}

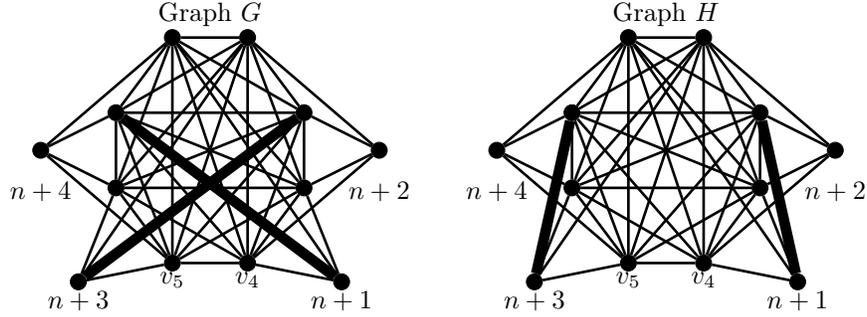
\begin{figure}[ht]
    \begin{multicols}{2}
    \centering
    Graph $G$
        \begin{tikzpicture}
            \tikzstyle{every node} = [inner sep = 0.0, minimum size = 6pt, draw = black]
            \node (v1) [color = white,fill = black, circle , draw = black] at (3.00,4.00) {};
            \node (v2) [color = white,fill = black, circle , draw = black] at (3.75,3.00) {};
            \node (v3) [color = white,fill = black, circle , draw = black] at (3.75,2.00) {};
            \node (v4) [color = white,fill = black, circle , draw = black] at (3.00,1.00) [label=below:$v_4$] {};
            \node (v5) [color = white,fill = black, circle , draw = black] at (2.00,1.00) [label=below:$v_5$] {};
            \node (v6) [color = white,fill = black, circle , draw = black] at (1.25,2.00) {};
            \node (v7) [color = white,fill = black, circle , draw = black] at (1.25,3.00) {};
            \node (v8) [color = white,fill = black, circle , draw = black] at (2.00,4.00) {};
            \node (v9) [color = white,fill = black, circle , draw = black] at (4.25,0.75) [label=below:$n+1$] {};
            \node (v10) [color = white,fill = black, circle , draw = black] at (4.75,2.50) [label={[label distance=3mm]below:$n+2$}] {};
            \node (v11) [color = white,fill = black, circle , draw = black] at (0.75,0.75) [label=below:$n+3$] {};
            \node (v12) [color = white,fill = black, circle , draw = black] at (0.25,2.50) [label={[label distance=3mm]below:$n+4$}] {};
            \draw [black,line width = 1]  (v1) -- (v2);
            \draw [black,line width = 1]  (v1) -- (v3);
            \draw [black,line width = 1]  (v1) -- (v4);
            \draw [black,line width = 1]  (v1) -- (v5);
            \draw [black,line width = 1]  (v1) -- (v6);
            \draw [black,line width = 1]  (v1) -- (v7);
            \draw [black,line width = 1]  (v1) -- (v8);
            \draw [black,line width = 1]  (v2) -- (v3);
            \draw [black,line width = 1]  (v2) -- (v4);
            \draw [black,line width = 1]  (v2) -- (v5);
            \draw [black,line width = 1]  (v2) -- (v6);
            \draw [black,line width = 1]  (v2) -- (v7);
            \draw [black,line width = 1]  (v2) -- (v8);
            \draw [black,line width = 1]  (v3) -- (v4);
            \draw [black,line width = 1]  (v3) -- (v5);
            \draw [black,line width = 1]  (v3) -- (v6);
            \draw [black,line width = 1]  (v3) -- (v7);
            \draw [black,line width = 1]  (v3) -- (v8);
            \draw [black,line width = 1]  (v4) -- (v5);
            \draw [black,line width = 1]  (v4) -- (v6);
            \draw [black,line width = 1]  (v4) -- (v7);
            \draw [black,line width = 1]  (v4) -- (v8);
            \draw [black,line width = 1]  (v5) -- (v6);
            \draw [black,line width = 1]  (v5) -- (v7);
            \draw [black,line width = 1]  (v5) -- (v8);
            \draw [black,line width = 1]  (v6) -- (v7);
            \draw [black,line width = 1]  (v6) -- (v8);
            \draw [black,line width = 1]  (v7) -- (v8);
            \draw [black,line width = 1]  (v10) -- (v1);
            \draw [black,line width = 1]  (v10) -- (v2);
            \draw [black,line width = 1]  (v10) -- (v3);
            \draw [black,line width = 1]  (v10) -- (v4);
            \draw [black,line width = 1]  (v9) -- (v8);
            \draw [black,line width = 4]  (v9) -- (v7);
            \draw [black,line width = 1]  (v9) -- (v3);
            \draw [black,line width = 1]  (v9) -- (v4);
            \draw [black,line width = 1]  (v12) -- (v8);
            \draw [black,line width = 1]  (v12) -- (v7);
            \draw [black,line width = 1]  (v12) -- (v6);
            \draw [black,line width = 1]  (v12) -- (v5);
            \draw [black,line width = 1]  (v11) -- (v1);
            \draw [black,line width = 4]  (v11) -- (v2);
            \draw [black,line width = 1]  (v11) -- (v6);
            \draw [black,line width = 1]  (v11) -- (v5);
        \end{tikzpicture}

    \centering
    Graph $H$
        \begin{tikzpicture}
        \tikzstyle{every node} = [inner sep = 0.0, minimum size = 6pt, draw = black]
            \node (v1) [color = white,fill = black, circle , draw = black] at (13.00,4.00) {};
            \node (v2) [color = white,fill = black, circle , draw = black] at (13.75,3.00) {};
            \node (v3) [color = white,fill = black, circle , draw = black] at (13.75,2.00) {};
            \node (v4) [color = white,fill = black, circle , draw = black] at (13.00,1.00) [label=below:$v_4$] {};
            \node (v5) [color = white,fill = black, circle , draw = black] at (12.00,1.00) [label=below:$v_5$] {};
            \node (v6) [color = white,fill = black, circle , draw = black] at (11.25,2.00) {};
            \node (v7) [color = white,fill = black, circle , draw = black] at (11.25,3.00) {};
            \node (v8) [color = white,fill = black, circle , draw = black] at (12.00,4.00) {};
            \node (v9) [color = white,fill = black, circle , draw = black] at (14.25,0.75) [label=below:$n+1$] {};
            \node (v10) [color = white,fill = black, circle , draw = black] at (14.75,2.50) [label={[label distance=3mm]below:$n+2$}] {};
            \node (v11) [color = white,fill = black, circle , draw = black] at (10.75,0.75) [label=below:$n+3$] {};
            \node (v12) [color = white,fill = black, circle , draw = black] at (10.25,2.50) [label={[label distance=3mm]below:$n+4$}] {};
            \draw [black,line width = 1]  (v1) -- (v2);
            \draw [black,line width = 1]  (v1) -- (v3);
            \draw [black,line width = 1]  (v1) -- (v4);
            \draw [black,line width = 1]  (v1) -- (v5);
            \draw [black,line width = 1]  (v1) -- (v6);
            \draw [black,line width = 1]  (v1) -- (v7);
            \draw [black,line width = 1]  (v1) -- (v8);
            \draw [black,line width = 1]  (v2) -- (v3);
            \draw [black,line width = 1]  (v2) -- (v4);
            \draw [black,line width = 1]  (v2) -- (v5);
            \draw [black,line width = 1]  (v2) -- (v6);
            \draw [black,line width = 1]  (v2) -- (v7);
            \draw [black,line width = 1]  (v2) -- (v8);
            \draw [black,line width = 1]  (v3) -- (v4);
            \draw [black,line width = 1]  (v3) -- (v5);
            \draw [black,line width = 1]  (v3) -- (v6);
            \draw [black,line width = 1]  (v3) -- (v7);
            \draw [black,line width = 1]  (v3) -- (v8);
            \draw [black,line width = 1]  (v4) -- (v5);
            \draw [black,line width = 1]  (v4) -- (v6);
            \draw [black,line width = 1]  (v4) -- (v7);
            \draw [black,line width = 1]  (v4) -- (v8);
            \draw [black,line width = 1]  (v5) -- (v6);
            \draw [black,line width = 1]  (v5) -- (v7);
            \draw [black,line width = 1]  (v5) -- (v8);
            \draw [black,line width = 1]  (v6) -- (v7);
            \draw [black,line width = 1]  (v6) -- (v8);
            \draw [black,line width = 1]  (v7) -- (v8);
            \draw [black,line width = 1]  (v10) -- (v1);
            \draw [black,line width = 1]  (v10) -- (v2);
            \draw [black,line width = 1]  (v10) -- (v3);
            \draw [black,line width = 1]  (v10) -- (v4);
            \draw [black,line width = 1]  (v9) -- (v8);
            \draw [black,line width = 4]  (v9) -- (v2);
            \draw [black,line width = 1]  (v9) -- (v3);
            \draw [black,line width = 1]  (v9) -- (v4);
            \draw [black,line width = 1]  (v12) -- (v8);
            \draw [black,line width = 1]  (v12) -- (v7);
            \draw [black,line width = 1]  (v12) -- (v6);
            \draw [black,line width = 1]  (v12) -- (v5);
            \draw [black,line width = 1]  (v11) -- (v1);
            \draw [black,line width = 4]  (v11) -- (v7);
            \draw [black,line width = 1]  (v11) -- (v6);
            \draw [black,line width = 1]  (v11) -- (v5);
        \end{tikzpicture}
    \end{multicols}
\caption{Graphs $G$ and $H$}
\label{fig:Graphs $G$ and $H$}
\end{figure}

\pf{We present an example depicting two infinite classes of non-isomorphic graphs with the same neighbor matrix below. 

We define graphs $G$ and $H$ in accordance with Table~\ref{table:noniso}}. Vertices $1$ to $n$ define a $K_n$. In Figure~\ref{fig:Graphs $G$ and $H$}, $n=8$.

    \begin{table}[ht]
    \centering 
        \begin{tabular}{c c c} 
            & $G$ & $H$ \\ \hline
            vertex $n+1$ adjacent to & $n,n-1,3,4$ & $n,2,3,4$  \\
            vertex $n+2$ adjacent to & $1,2,3,4$ & $1,2,3,4$ \\
            vertex $n+3$ adjacent to & $n,n-1,n-2,n-3$ & $n,n-1,n-2,n-3$ \\
            vertex $n+4$ adjacent to & $1,2,n-2,n-3$ & $1,n-1,n-2,n-3$ \\ \hline
        \end{tabular}
    \caption{Non-isomorphic Graphs $G$ and $H$ with Identical Neighbor Matrices} 
    \label{table:noniso} 
    \end{table}

Figure~\ref{fig:Graphs $G$ and $H$} is an illustration of the graphs constructed in Table~\ref{table:noniso}. The edges drawn in heavier pitch are those edges which are different between graphs $G$ and $H$.


Given $G$ and $H$, $X^{dist-k}(G) = X^{dist-k}(H).$ In graph $G$, each degree-$4$ vertex has a maximum of two dist-$1$ neighbors shared with any other degree-$4$ vertex. In graph $H$, each degree-$4$ vertex has one vertex with which it shares three dist-$1$ neighbors. Therefore $G \centernot\cong H.$~\qed

It is clear that the implication stated in Proposition~\ref{Prop:Congruent} holds for trees. Buckley and Harary provide a counterexample for Remark~\ref{rmk:different} which may be extended to show identical neighbor matrices do not imply isomorphism for trees~\cite{buckley1990distance}. As algorithms exist to construct neighbor matrices, the contrapositive of Theorem~\ref{Prop:Congruent} is a method to verify the non-isomorphism of two graphs.

\begin{figure}[ht]
    \centering
    \begin{multicols}{2}
        \begin{tikzpicture}[scale=0.2]
            \tikzstyle{every node} = [inner sep = 0.0, minimum size = 6pt, draw = black]
            \node (v1) [color = white,fill = black, circle , draw = black] at (10.00,20.00) [label= {[label distance = 2mm]$v_{1}$}] {};
            \node (v2) [color = white,fill = black, circle , draw = black] at (20.0,10.00) [label= {[label distance = 2mm]right:$v_{2}$}] {};
            \node (v3) [color = white,fill = black, circle , draw = black] at (15.0,0) [label= {[label distance = 2mm]below:$v_{3}$}] {};
            \node (v4) [color = white,fill = black, circle , draw = black] at (5.00,0) [label= {[label distance = 2mm]below:$v_{4}$}] {};
            \node (v5) [color = white,fill = black, circle , draw = black] at (0,10.00) [label= {[label distance = 2mm]left:$v_{5}$}] {};
            \node (u1) [color = white,fill = black, circle , draw = black] at (10.00,15.00) [label= {[label distance = 2mm]left:$u_{1}$}] {};
            \node (u2) [color = white,fill = black, circle , draw = black] at (15.0,10.00) [label= {[label distance = 2mm]above:$u_{2}$}] {};
            \node (u3) [color = white,fill = black, circle , draw = black] at (12.50,3.5) [label= {[label distance = 2mm]below:$u_{3}$}] {};
            \node (u4) [color = white,fill = black, circle , draw = black] at (7.50,3.5) [label= {[label distance = 2mm]below:$u_{4}$}] {};
            \node (u5) [color = white,fill = black, circle , draw = black] at (5,10) [label= {[label distance = 2mm]above:$u_{5}$}] {};
            \draw [dashed]  (u1) -- (u3);
            \draw [black,line width = .2]  (u3) -- (u5);
            \draw [black,line width = .2]  (u5) -- (u2);
            \draw [black,line width = .2]  (u2) -- (u4);
            \draw [dashed]  (u4) -- (u1);
            \draw [dotted, ultra thick]  (v1) -- (v2);
            \draw [dashed]  (v2) -- (v3);
            \draw [black,line width = .2]  (v3) -- (v4);
            \draw [dashed]  (v4) -- (v5);
            \draw [dotted, ultra thick]  (v5) -- (v1);
            \draw [dotted, ultra thick]  (u1) -- (v1);
            \draw [dashed]  (u2) -- (v2);
            \draw [black,line width = .2]  (u3) -- (v3);
            \draw [black,line width = .2]  (u4) -- (v4);
            \draw [dashed]  (u5) -- (v5);
        \end{tikzpicture}
    \columnbreak
        \begin{tikzpicture}[scale=0.2]
            \tikzstyle{every node} = [inner sep = 0.0, minimum size = 6pt, draw = black]
            \node (v1) [color = white,fill = black, circle , draw = black] at (10.00,20.00) [label= {[label distance = 2mm]$v_{1}$}] {};
            \node (v2) [color = white,fill = black, circle , draw = black] at (20.0,10.00) [label= {[label distance = 2mm]right:$v_{2}$}] {};
            \node (v3) [color = white,fill = black, circle , draw = black] at (15.0,0) [label= {[label distance = 2mm]below:$v_{3}$}] {};
            \node (v4) [color = white,fill = black, circle , draw = black] at (5.00,0) [label= {[label distance = 2mm]below:$v_{4}$}] {};
            \node (v5) [color = white,fill = black, circle , draw = black] at (0,10.00) [label= {[label distance = 2mm]left:$v_{5}$}] {};
            \node (u1) [color = white,fill = black, circle , draw = black] at (10.00,15.00) [label= {[label distance = 2mm]left:$u_{1}$}] {};
            \node (u2) [color = white,fill = black, circle , draw = black] at (15.0,10.00) [label= {[label distance = 2mm]above:$u_{2}$}] {};
            \node (u3) [color = white,fill = black, circle , draw = black] at (12.50,3.5) [label= {[label distance = 2mm]below:$u_{3}$}] {};
            \node (u4) [color = white,fill = black, circle , draw = black] at (7.50,3.5) [label= {[label distance = 2mm]below:$u_{4}$}] {};
            \node (u5) [color = white,fill = black, circle , draw = black] at (5,10) [label= {[label distance = 2mm]above:$u_{5}$}] {};
            \draw [black,line width = 0.2]  (u1) -- (u3);
            \draw [black,line width = 0.2]  (u3) -- (u5);
            \draw [dashed]  (u5) -- (u2);
            \draw [dashed]  (u2) -- (u4);
            \draw [black,line width = 0.2]  (u4) -- (u1);
            \draw [dotted, ultra thick]  (v1) -- (v2);
            \draw [dotted, ultra thick]  (v2) -- (v3);
            \draw [dashed]  (v3) -- (v4);
            \draw [black,line width = 0.2]  (v4) -- (v5);
            \draw [dashed]  (v5) -- (v1);
            \draw [dashed]  (u1) -- (v1);
            \draw [dotted, ultra thick]  (u2) -- (v2);
            \draw [dashed]  (u3) -- (v3);
            \draw [black,line width = 0.2]  (u4) -- (v4);
            \draw [black,line width = 0.2]  (u5) -- (v5);
        \end{tikzpicture}
    \end{multicols}
\caption{Two views of distance-$1$ and distance-$2$ neighbors (from $v_1$ and $v_2$) in the Petersen graph} \label{fig:Petersen}
\end{figure}
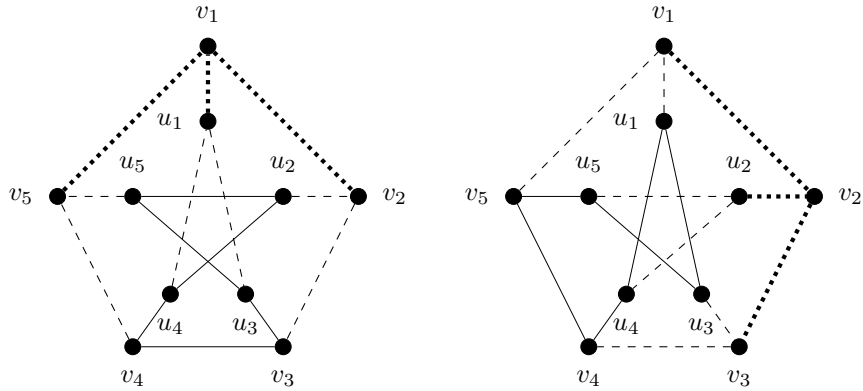

The beauty of the newly introduced neighbor matrix extends from how naturally it captures the topology of an arbitrary graph, coupled with the simplicity of calculating matrix entries. A neighbor matrix of a graph captures each vertex's ``view'' of the graph through shortest paths, thus reaching every vertex of the graph, implicitly accounting for edges, cycles and other subgraphs of $G$ (see Figure~\ref{fig:Petersen}, two views of the Petersen graph).

From the topological information encoded in the neighbor matrix, we may extract many graph parameters, as shown in Theorem~\ref{prop:properties}.

\begin{theorem}\label{prop:properties}
For a simple graph $G$ the following topological information can be observed from the neighbor matrix $X^{dist-k}(G)=[x_{ij}], (1 \le i \le n, 1 \le j\le k)$:
    \begin{enumerate}[label=(\alph*), labelindent=\parindent, leftmargin=*, start=1]
        \item $G$ is connected if and only if $\displaystyle{\sum_{i=1}^n\sum_{j=1}^k x_{ij} = n(n-1)}.$
            \begin{itemize}
                \item The number of components in $G$, $\kappa(G)$, is $\displaystyle{\sum_{d=1}^{n-1} \frac{|\mathscr{S}_d|}{d+1}}, $ where\\ $\displaystyle{\mathscr{S}_d=\{v_i:\sum_{j=1}^k v_{ij} = d\}}.$
            \end{itemize}
        \item The radius of the graph is $Rad(G) = \max\{j:  x_{ij} \neq 0, \forall i , 1 \le i \le n\}.$
        \item The center of the graph is $Cen(G) = G[\{v_{i}: x_{i\rad(G)} \neq 0, x_{ij} = 0, \forall j>rad(G)\}].$
        \item The closeness centrality of vertex $i$ in a connected graph is $$\displaystyle{CC_i = \frac{n-1}{\sum_{j=1}^k j \cdot x_{ij}}}.$$ 
        \item The average distance between vertices in a connected graph $G$ is $$\displaystyle {\frac{1}{n(n-1)} \sum_{i=1}^n\sum_{j=1}^k (j\cdot x_{ij}) = \frac{1}{n(n-1)} \sum_{j=1}^k\sum_{i=1}^n (j\cdot x_{ij})}.$$
        \item  Given $\diam(G) = k$ (the number of columns in $X^{dist-k}(G)),$ $ $ $(1 \le k \le n-1)$, the graph periphery is $Per(G) = G[\{v_{i}: x_{ik} \neq 0\,\ \ \forall i,\ \ 1 \leq i \leq n\}].$
        \item  The first column is a representation of the degree sequence, in which the degree centrality of vertex $i$ is the first entry in row $i$, $x_{i1}$. The number of edges in the graph is the half-sum of the entries in the first column:  $$\displaystyle{m=\frac{1}{2} \cdot \sum_{i=1}^n x_{i1}}.$$
            \begin{itemize}
            \item The \textit{density} of a graph is the ratio of possible edges to the edges that are actually present in the graph: 
    $${\frac{m}{\binom{n}{2}}  = \frac{2m}{n(n-1)}}.$$
            \end{itemize}
        \item For each column $j \in X^{dist-k}(G)$, $\displaystyle{\sum_{i=1}^n x_{ij} = 2\left(|E(G^j)| - |E(G^{j-1})|\right)}.$ 
        \item For each $s$ $\big(1 \le s \le \diam(G) \big)$, the number of edges in the power graph $G^s$ is given by $$\sum_{j=1}^s \sum_{i=1}^n x_{ij} = 2|E(G^s)|.$$ 
    \end{enumerate}
\end{theorem}

\pf{
\begin{itemize}
    \item [(a)]  Suppose $G$ is connected. Fix a vertex $v$, and note that $\forall u \in V(G-v),$ a geodesic exists such that $d(u, v) = j $ $ (1 \le j \le k)$. There are $n-1$ choices for $u$ in $G-v$. By construction of $X^{dist-k}(G)$, the length of a geodesic between all vertex pairs $u, v \in V(G)$ is counted. So, for each of the $n$ rows of the matrix there are $n-1$ distance entries, totalling $n(n-1)$. For the converse, assume, to the contrary, that $G$ is not connected and $\displaystyle{\sum_{i=1}^n\sum_{j=1}^k x_{ij} = n(n-1)}$. Then, there is a component $G'$ containing $q$ vertices $(1 \leq q \leq n-1)$; we choose $v \notin V(G')$. There are at most $n-1-q$ vertices adjacent to $v$. By construction of $X^{dist-k}(G)$ the entries of the row corresponding to vertex $v$ sum to $n-1-q$, and also there are at most $n-2$ geodesics that are counted in the other rows  of $X^{dist-k}(G)$.  Therefore
$\displaystyle{\sum_{i=1}^n\sum_{j=1}^k x_{ij} \leq}$ $\displaystyle{(n-1)(n-2) + (n-1-q)  < n(n-1)}$, which contradicts the initial statement that $\displaystyle{\sum_{i=1}^n\sum_{j=1}^k x_{ij} = n(n-1)}$.~\qed{}

    \begin{itemize}
    \item We proved in (a) that a connected graph (i.e., graph with $1$ component) has $n$ neighbor matrix row sums of $n-1$. It follows that each connected component of order $d$ will have $d$ neighbor matrix row sums of $d-1$. Therefore, given $\mathscr{S}_d=\{v_i:\sum_{j=1}^k v_{ij} = d\}$, $\forall d: 1 \leq d \leq n-1$, the number of connected components in $G$ of order $d$ is $$\displaystyle{\kappa(G)_d = \frac{|\mathscr{S}_d|}{d+1}},$$ and, $\kappa(G) = \displaystyle{\sum_{d=1}^{n-1} \kappa(G)_d}$.~\qed{}
    \end{itemize}

    By construction of $X^{dist-k}$, each value $x_{ij}$ is the number of geodesics originating at vertex $i$ of length $j \, (1 \leq j \leq k)$ with different terminal vertices. We weight the entries of the $k$-matrix row by multiplying each value $x_{ij}$ by its associated distance, $j$. As seen in (a), the number of possible geodesics in $G$ in a simple, connected graph is $n(n-1)$. 
    \begin {itemize}
        \item [(d)] After summing over row $i$, division by the number of geodesics provides the weighted average which is the average distance in $G$ from vertex $i$ to all other vertices. As stated in Section~\ref{section:definitions}, the closeness centrality is the inverse of this average distance.~\qed{}
        \item [(e)] After summing over rows and columns, division by the number of geodesics provides the weighted average which is the average distance between vertices in $G$.~\qed{}
    \end{itemize}
    \item [(h)] The left hand side of the equation double counts the sum of distance $j$ neighbors ($1 \le j \le \diam(G)$) from each vertex in $G$.  The right hand side does the same by calculating the number of vertex pairs of distance at most $j$ apart, and then it removes the number of vertex pairs of distance at most $j-1$. \qed{}
    \item [(i)] The left hand side of the equation double counts the sum of distance 1, 2, \ldots, $j$ neighbors ($1 \le j \le \diam(G)$) from each vertex in $G$.  The right hand side does the same by calculating the number of vertex pairs of distance at most $j$ apart.~\qed{}
\end{itemize}
}

\section{The Neighbor Matrix and Graph Orbits}\label{sec:orbits}

As stated in Section~\ref{section:definitions}, vertices in the same orbit of a graph form an equivalence class. These equivalence classes are related to neighbor matrix rows, as stated in Proposition~\ref{prop:orbitRow}.

\begin{proposition}\label{prop:orbitRow}
Given vertices $v_i$ and $v_j$ in graph $G$ we have $$o_G(i) = o_G(j) \implies X_i^{dist-k}(G) = X_j^{dist-k}(G).$$

\end{proposition}

\pf{Suppose otherwise, that $v_i$ and $v_j$ are in the same orbit and $X_i^{dist-k}(G)$ $\neq X_j^{dist-k}(G)$. Then, there is at least one distance value {\it f} for which $v_i$ has more (or less) $distance-f$ neighbors than $v_j$. Therefore, there is no automorphism that maps $v_i$ to $v_j$ in $Aut(G)$. This is a contradiction, as $v_i$ and $v_j$ are in the same orbit, such a permutation of vertices must indeed exist in $Aut(G)$.\qed{}

\begin{remark}
The converse of Proposition~\ref{prop:orbitRow} is false. That is, $$X_i^{dist-k}(G) = X_j^{dist-k}(G) \centernot \implies o_G(i) = o_G(j).$$
\end{remark}

\pf{The graph $G$ in Figure \ref{fig:Orbits} shows an example of a graph with a cycle. Note that vertices $v_i$ and $v_j$ have identical neighbor matrix rows. There are two different paths from vertex $v_i$ to vertex $v_k$, a distance-$2$ neighbor. All paths from vertex $v_j$ to distance-$2$ neighbors are unique; therefore there is no automorphism in $Aut(G1)$ that maps $v_i$ to $v_j$ and the two vertices are in different orbits.  

Consider graph $T$ in Figure~\ref{fig:Orbits}. Note that vertices $v_y$ and $v_z$ have identical neighbor matrix rows.  The subgraphs rooted at vertices $v_y$ and $v_z$, respectively each have two distance-$2$ neighbors. In the subgraph rooted at $v_y$, the paths to the two distance-$2$ neighbors share a common vertex. In the subgraph rooted at $v_z$, the paths to the two distance-$2$ neighbors do not share a common vertex. Therefore, there is no automorphism in $Aut(T)$ that maps $v_y$ to $v_z$ and the two vertices are in different orbits.\qed{} 

}
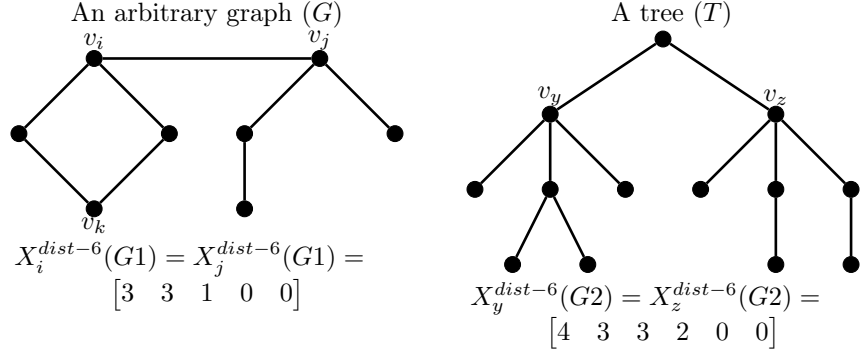
\begin{figure}[ht]
    \begin{multicols}{2}
    \centering
        An arbitrary graph ($G$)
            \begin{tikzpicture}
                \tikzstyle{every node} = [inner sep = 0.0, minimum size = 6pt, draw = black]
                \node (v1) [color = white,fill = black, circle , draw = black] at (0.00,0.00) [label=above:$v_i$] {};
                \node (v2) [color = white,fill = black, circle , draw = black] at (3.0,0.00) [label=above:$v_j$] {};
                \node (v3) [color = white,fill = black, circle , draw = black] at (-1.0,-1.0) {};
                \node (v4) [color = white,fill = black, circle , draw = black] at (1.0,-1.0) {};
                \node (v5) [color = white,fill = black, circle , draw = black] at (2.00,-1.0) {};
                \node (v6) [color = white,fill = black, circle , draw = black] at (4.0,-1.0) {};
                \node (v7) [color = white,fill = black, circle , draw = black] at (0,-2.00) [label=below:$v_k$] {};
                \node (v8) [color = white,fill = black, circle , draw = black] at (2.00,-2.00) {};
                \draw [black,line width = 1]  (v1) -- (v2);
                \draw [black,line width = 1]  (v1) -- (v3);
                \draw [black,line width = 1]  (v1) -- (v4);
                \draw [black,line width = 1]  (v2) -- (v5);
                \draw [black,line width = 1]  (v2) -- (v6);
                \draw [black,line width = 1]  (v3) -- (v7);
                \draw [black,line width = 1]  (v4) -- (v7);
                \draw [black,line width = 1]  (v5) -- (v8);
            \end{tikzpicture}
                    
    $X_i^{dist-6}(G1) = X_j^{dist-6}(G1) =$
    \newline
    $\begin{bmatrix}
        3 & 3 & 1 & 0 & 0
    \end{bmatrix}$
                    
    \columnbreak
            
    \centering
    \     \ A tree ($T$)
        \begin{tikzpicture}
            \tikzstyle{every node} = [inner sep = 0.0, minimum size = 6pt, draw = black]
            \node (u1) [color = white,fill = black, circle , draw = black] at (-0.50,0.00) {};
            \node (u2) [color = white,fill = black, circle , draw = black] at (-2.0,-1.00) [label=above:$v_y$] {};
            \node (u3) [color = white,fill = black, circle , draw = black] at (1.00,-1.00) [label=above:$v_z$] {};
            \node (u4) [color = white,fill = black, circle , draw = black] at (-3.00,-2.00) {};
            \node (u5) [color = white,fill = black, circle , draw = black] at (-2.00,-2.00) {};
            \node (u6) [color = white,fill = black, circle , draw = black] at (-1.0,-2.00) {};
            \node (u7) [color = white,fill = black, circle , draw = black] at (0.0,-2.00) {};
            \node (u8) [color = white,fill = black, circle , draw = black] at (1.00,-2.00) {};
            \node (u9) [color = white,fill = black, circle , draw = black] at (2.00,-2.00) {};
            \node (u10) [color = white,fill = black, circle , draw = black] at (-2.50,-3.00) {};
            \node (u11) [color = white,fill = black, circle , draw = black] at (-1.50,-3.00) {};
            \node (u12) [color = white,fill = black, circle , draw = black] at (1.00,-3.00) {};
            \node (u13) [color = white,fill = black, circle , draw = black] at (2.00,-3.00) {};
            \draw [black,line width = 1]  (u1) -- (u2);
            \draw [black,line width = 1]  (u1) -- (u3);
            \draw [black,line width = 1]  (u2) -- (u4);
            \draw [black,line width = 1]  (u2) -- (u5);
            \draw [black,line width = 1]  (u2) -- (u6);
            \draw [black,line width = 1]  (u3) -- (u7);
            \draw [black,line width = 1]  (u3) -- (u8);
            \draw [black,line width = 1]  (u3) -- (u9);
            \draw [black,line width = 1]  (u5) -- (u10);
            \draw [black,line width = 1]  (u5) -- (u11);
            \draw [black,line width = 1]  (u8) -- (u12);
            \draw [black,line width = 1]  (u9) -- (u13);
        \end{tikzpicture}
                                    
    $X_y^{dist-6}(G2) = X_z^{dist-6}(G2) =$
    \newline
    $\begin{bmatrix}
        4 & 3 & 3 & 2 & 0 & 0
    \end{bmatrix}$
                    
    \end{multicols}
\caption{Graphs $G1$ and $G2$, Counterexamples to Proposition~\ref{prop:orbitRow}}
\label{fig:Orbits}
\end{figure}

\section{Finding the Neighbor Matrix}\label{section:k-matrix}

The neighbor matrix may be fully determined using powers of adjacency matrices, or any of a variety of computer algorithms that determine all-pairs shortest paths. Propositions~\ref{Prop:graph_powers} and~\ref{Prop:matrix1} present two algebraic methods of obtaining a neighbor matrix representation of $G$.  For these propositions we use $X_{A}^{dist-k}$ to refer to the neighbor matrix obtained from an adjacency matrix $A$ of $G$, and $\vec{1}$ to represent the $n \times 1$ vector whose entries are $1$.
\bigskip

\begin{proposition}\label{Prop:graph_powers}
Each neighbor matrix can be obtained using adjacency matrices of powers of the original graph $G$, as shown in equation~$\ref{eqn:powers}$. 
Entries of  $-1$ in $A(G^i) - A(G^{i-1})$ are replaced by $0$ before multiplication by $\vec{1}$:

    \begin{equation}\label{eqn:powers}
        X_{A}^{dist-k} =  
        \begin{bmatrix}
            A(G)\vec{1} & \left(A(G^2)-A(G)\right) \vec{1}  & \ldots &  \left(A(G^k) - A(G^{k-1})\right) \vec{1} $ $ \\
        \end{bmatrix}.
    \end{equation}
\end{proposition}

\pf
Notice that from the definition of $G^k$, the adjacency matrix of $G^k$, namely $A(G^k)$, includes entries of $1$ for all pairs of vertices of distance $k$ or less in the original graph $G$. Therefore, when we subtract $A(G^{k-1})$ from $A(G^k)$ and replace the $-1$ entries by $0$, we are left with a $(0,$ $1)$-matrix with $1$ entries representing pairs of vertices exactly distance $k$ apart in $G$. Multiplication by $\vec{1}$ creates a vector consisting of the number of 
vertices at distance $k$ from each vertex in $G$.~\qed

\bigskip

We let $\mathscr{A}^j$ be the Boolean matrix obtained from
$\mathscr{A}^j = (A^j)^b  - \mathscr{A}^{j-1}$, where $ j \ge 2$, by replacing the value of $-1$ and diagonal entries with zero; $\mathscr{A}^1 = A$, the adjacency matrix.

\bigskip 
\begin{proposition}
\label{Prop:matrix1}
    Each unsorted neighbor matrix, $X^{dist-k(u)}$ can be obtained using the boolean matrices $\mathscr{A}^j (1 \le j \le k )$, as shown in equation~$\ref{eqn:boolean}$
    \begin{equation}\label{eqn:boolean}
        X_{\mathscr{A}}^{dist-k(u)}(G) =  
        \begin{bmatrix}
            \mathscr{A}^1\vec{1} & \mathscr{A}^2 \vec{1} &  \mathscr{A}^3 \vec{1}   & \ldots &\mathscr{A}^k \vec{1}\\
        \end{bmatrix}.
    \end{equation}
\end{proposition}

\pf
Notice that multiplying the adjacency matrix $\mathscr{A}^1$ of $G$ by the column vector $\vec{1}$ we obtain the first column of the matrix $X^{dist-k}(G)$, the degree sequence. The second column of $X^{dist-k}(G)$ is $\mathscr{A}^2  \vec{1} - \mathscr{A} \vec{1}$, as it counts all the vertices of distance $2$ or less, and then it subtracts the vertices of distance $1$, i.e. the neighbors of each fixed vertex in $V(G)$, as shown in Theorem~\ref{prop:properties}.  Similarly, $\mathscr{A}^k \vec{1} - \mathscr{A}^{k-1} \vec{1}$ counts the number of vertices $k$-hops away, and it subtracts the number of vertices $k-1$ hops away.~\qed

After a reverse lexicographic sort, $X^{dist-k(u)}$ becomes the neighbor matrix $X^{dist-k}$.

\bigskip

We next show that the neighbor-matrices $X^{dist-k}$ obtained from different adjacency matrices of a graph $G$ are related through the same permutation matrices as the adjacency matrices themselves.  

\bigskip
\begin{proposition}
\label{Prop:matrix2}
    Given $A$ and $B$, two different adjacency matrices of graph $G$ (so $B = P \cdot A \cdot P^T$, for some permutation matrix $P$),   $$X_{{B}}^{dist-k}(G) = P \cdot X_{{A}}^{dist-k}(G) \cdot P^T.$$
\end{proposition}

\pf  
Let  $A$ and $B$ be two different adjacency matrices of $G$, so $B = P \cdot A \cdot P^T$ for some permutation matrix $P$.  Since $P$ is an orthogonal matrix, we have that
$B^i = (P A P^T)^i = P A^i P^T$, for $1 \le i \le k$.  Therefore $B^j-B^{j-1} = P (B^j-A^{j-1}) P^T$  for each $j (1 \le j \le k)$, which implies that $X_B^{dist-k}(G) = P \cdot X_A^{dist-k}(G) \cdot P^T.$~\qed

\section{Applications}\label{section:applications}
There are many potential applications for the neighbor matrix, as it is a compact structure and rich in topological information. Two such applications are graph comparison and the identification of vertices with high topological significance, relative to other vertices in the graph.

\subsection{Graph Comparison}\label{subsection:comparison}
The neighbor matrix enables the analysis of two graphs and insight into their topological {\it similarity}, which we consider to be a relaxation of isomorphism. We consider two graphs to be {\it similar} if they have the same distance distributions, i.e., their neighbor matrices are the same. For example, using this definition, graphs $G$ and $H$ in Figure~\ref{fig:Graphs $G$ and $H$} are {\it similar}.

It has previously been established that using individual characteristics provides an incomplete picture of a graph's structure \cite{smetric1, smetric2}. This is because information is missing that is needed to fully characterize the graph. Several measures follow that illustrate the problem of different structures with the same (or close) comparison criteria:

\begin{itemize}
    \item Degree Sequence: a cycle $C_{3k}$ and the disconnected graph obtained from the  union of $k$ copies of the $C_3$ cycle (or, $k C_3$) each have constant, identical degree sequences
    \item Average Distance:  a wheel on $32$ vertices $W_{1,31}$ and the $16$-vertex star graph $K_{1, 15}$ have the average  path length of $1.875$. 
    \item Average Clustering Coefficient (close, positive): a wheel on $21$ vertices $W_{1, 20}$ has the average clustering coefficient of $0.63888$, while the $(7,7)$ Barbell graph (two copies of $K_7$ connected by a path with $7$ additional vertices) has average clustering coefficient of $.63945$.
    \item Average Clustering Coefficient (same): a graph of small diameter (star) and a large diameter (path) each have a clustering coefficient of $0$.
    \item Pearson Coefficient:  a wheel on $n+1$ vertices ($n \ge 5$) $W_{1,n}$ and the $(3,5)$ Barbell graph (two copies of $K_3$ connected by a path with  $5$ more vertices) have the Pearson correlation coefficient of $-.33333$. 
    \item $s$-metric: for each $r-$regular graph, the normalized $s$-metric is $1$ because there is no variability in the degree sequence in the background set of graphs.

\end{itemize}

Though variability of graphs with the same or similar graph invariant values has been observed before, it has not been adequately addressed. A widespread desire remains in the community of decision makers operating in a connected world to compare graph and network structures in a scientifically-viable manner. This, coupled with the well-documented lack of a comprehensive, intuitive metric useful in making this comparison served as a motivation for development of the neighbor matrix to describe graphs.

Consider Table~\ref{table:Measures}, which includes several classes of graphs that are defined in Table~\ref{table:graphs}. The entries in bold are the graph invariants previously discussed that have identical values when calculated from different graphs. Note that eight of the graphs represented in Table~\ref{table:Measures} have the same average clustering coefficient: $0$; but neighbor matrix dimensions are different. This phenomenon persists through the other examples shown. A similar result comes from looking at the dimension of neighbor matrices associated with graphs manifesting the same average path length, Pearson correlation coefficent, and normalized s-metric. So, simple inspection seems to indicate utility in the neighbor matrix as a topological discriminator where other invariants fail, even before more detailed analysis.

In linear algebra, we use norms to determine the ``size'' of a matrix. As we see in Table~\ref{table:Measures}, the Frobenius norm is unique for each graph present, even those for which multiple measures are the same or even undefined. This uniqueness provides the foundation for future work that uses the norm to calculate a graph comparison metric that is potentially more topologically insightful than measures currently in use.

\begin{table}[ht]
\centering
    \begin{tabular}{p{0.75cm} p{4.5cm} p{0.75cm} p{4.5cm}} \hline
        Graph & Description & Graph & Description \\\hline \hline
        $B_{7,7}$ & 2 x $K_7$ connected by $P_7$ & $W_{1,20}$ & wheel on $21$ vertices\\\hline
        $K_{1,15}$ & star on $16$ vertices & $P_{16}$ & path on $16$ vertices \\\hline
        $W_{1,31}$ & wheel on $32$ vertices & $B_{10,12}$ & 2 x $K_{10}$ connected by $P_{12}$ \\\hline
        $K_{32}$ & complete graph on $32$ vertices & $K_{16,16}$ & $32$-vertex complete bipartite \\\hline
        $CL_{32}$ & circular ladder on 32 vertices & $C_{32}$ & cycle on $32$ vertices \\\hline
        $H_{32}$ & hypercube on $32$ vertices & $L_{21,11}$ & $K_{21}$ connected to $P_{11}$ \\\hline
        $P_{32}$ & path on $32$ vertices & $K_{1,31}$ & star on $32$ vertices \\\hline \hline
    \end{tabular}
\caption{Graphs Analyzed} 
\label{table:graphs} 
\end{table}

\begin{table}[ht]
\centering
    \begin{tabular}{p{0.75cm} p{0.75cm} p{1.15cm} p{1.15cm} p{1.2cm} p{1.3cm} p{1.5cm}} \hline
        \small{Graph} & \small{AVG Distance} & \small{AVG Cluster Coeff} & \small{Pearson Correl Coeff} & \small{s-metric (Norm'd)} & \small{dimension ($X^{dist-k}$)} & \small{$X^{dist-k}$ Frobenius Norm} \\\hline \hline
        $B_{7,7}$ & 5.0 & {\textbf{0.640}} & 0.719 & 0.981 & 21 x 10 & 43.36 \\\hline
        $W_{1,20}$ & 1.81 & {\textbf{0.640}} & {\textbf{-0.333}} & 0.323 & 21 x 2 & 79.75 \\\hline
        $K_{1,15}$ & {\textbf{1.875}} & \textbf{0} & -1 & 0.133 & 16 x 2 & 56.39 \\\hline
        $P_{16}$ & 5.330 & \textbf{0} & -0.077 & 0.981 & 15 x 14 & 17.55 \\\hline
        $W_{1,31}$ & {\textbf{1.875}} & 0.648 & {\textbf{-0.333}} & 0.207 & 32 x 2 & 159.8 \\\hline
        $B_{10,12}$ & 7.323 & 0.613 & 0.866 & 0.990 & 32 x 15 & 77.05 \\\hline
        $K_{32}$ & 1.0 & 1.0 & * & {\textbf{1}} & 32 x 1 & 175.4 \\\hline
        $K_{16,16}$ & 1.484 & \textbf{0} & \textbf{*} & {\textbf{1}} & 32 x 2 & 124.1 \\\hline
        $CL_{32}$ & 4.645 & \textbf{0} & \textbf{*} & {\textbf{1}} & 32 x 9 & 60.66 \\\hline
        $C_{32}$ & 8.258 & \textbf{0} & \textbf{*} & {\textbf{1}} & 32 x 16 & 44.18 \\\hline
        $H_{32}$ & 2.581 & \textbf{0} & \textbf{*} & {\textbf{1}} & 32 x 5 & 89.62 \\\hline
        $L_{21,11}$ & 4.105 & 0.653 & 0.942 & 0.998 & 32 x 12 & 116.0 \\\hline
        $P_{32}$ & 11.0 & \textbf{0} & -0.033 & 0.992 & 32 x 31 & 38.37 \\\hline
        $K_{1,31}$ & 1.94 & \textbf{0} & -1 & 0.064 & 32 x 2 & 170.0 \\\hline \hline
    \end{tabular}
\caption{Graph Comparison Metrics} 
\flushleft{*The Pearson correlation coefficient is undefined in regular graphs because the denominator of the calculation is a variance of zero.}
\label{table:Measures} 
\end{table}

\clearpage

\subsection{Application 2: Identification of Topologically Influential Vertices}\label{subsection:influence}
Analysis using the topological information stored in the newly-developed neighbor matrix can rank the vertices in a graph based upon the structural importance of that vertex. If the removal of vertex $v_i$ results in a relatively large change to the distance distributions in a graph, we may assert $v_i$ is structurally important to the graph. Relatively small change to the distances in a graph upon vertex removal indicates less structural importance. There are other measures of topological importance currently in use---the centralities described in Section~\ref{section:definitions}; however, these measures do not directly address how vertices impact the distance between vertices in a graph.

This information can be leveraged towards recommendations for network attack and defense. In the attack, the attacker would choose to interdict the vertex or vertices that cause the greatest structural change. On the defense, the defender would choose to protect or ``harden'' those nodes that cause greatest structural change in the network.

\section{Conclusion}\label{section:conclusion}
This work introduced and explored the neighbor matrix as an algebraic structure that contains significant graph descriptive and topological information. We proved this topological richness by proving the presence of $11$ graph invariants in the neighbor matrix and relating the neighbor matrix to graph orbits. Though it does not inform isomorphism in arbitrary graphs, the neighbor matrix does provide a technique to verify the non-ismorphism of two graphs. The neighbor matrix has the potential to enable greater understanding of ``graph space,'' as it is simultaneously more compact and richer in information than current structures used in graph exploration. For example, the degree sequence has long been used as a foundational element for the exploration of families of graphs. Since the neighbor matrix extends the degree sequence through all the distances that comprise the graph, it promises to enhance our current capability to model, analyze, and understand graphs in all distance dimensions. Further exploration is warranted into how we may use the neighbor matrix to further refine models that use the degree sequence as an input, such as the configuration model of graph generation. We proved the neighbor matrix to contain eleven of the most commonly used invariants used to describe graphs. Future work will leverage this topological richness toward  generating insights into graph comparison and further, toward informing the identification of vertices whose influence on graph topology is significant. The neighbor matrix is a mathematically-manipulable structure that may move us toward a single characterization of certain classes of graphs, such as is found in the idea of a ``graph signature'' \cite{signature1,signature2}. Finally, as the neighbor matrix presents the sequences of all distances which manifest in a graph, it may provide a tool to verify the unigraphicality of certain classes of graphs.

\section{Acknowledgements}
This work was partially funded by a grant from the Department of Defense.


\bibliography{citations}

\end{document}